\documentclass[a4paper,11pt]{article}

\usepackage{amsmath}
\usepackage{graphicx}
\usepackage{amsfonts}
\usepackage{amssymb}
\usepackage{epsfig}
\usepackage{color}
\usepackage{psfrag}

 \newtheorem{theorem}{Theorem}%[section]
\newcommand{\EQ}{\begin{equation}}
\newcommand{\EN}{\end{equation}}

\newcommand{\bc}{\begin{center}}
\newcommand{\ec}{\end{center}}
\def\ba#1{\begin{array}{#1}\displaystyle}
\newcommand{\ea}{\end{array}}

\newcommand{\beq}{\begin{equation}}
\newcommand{\eeq}{\end{equation}}
\newcommand{\beqa}{\begin{eqnarray}}
\newcommand{\eeqa}{\end{eqnarray}}

\newcommand{\bi}{\begin{itemize}}
\newcommand{\ei}{\end{itemize}}
\def\t#1{\tilde{#1}}

\newcommand{\p}{\partial}

\newcommand{\Z}{{\mathbb{Z}}}

\newcommand{\R}{{\mathbb{R}}}

\newcommand{\hC}{{\hat{\mathbb{C}}}}

\newcommand{\onto}{\twoheadrightarrow}

\newcommand{\id}{{\rm id}}

\def\cl#1{\overline{#1}}

\begin{document}

\begin{center}
{\Large {\bf
Factorisation of conformal maps on finitely connected domains
}

\vspace{1cm}

Benjamin Doyon}

Department of Mathematics, King's College London\\
Strand, London, U.K.\\
benjamin.doyon@kcl.ac.uk

\end{center}

\vspace{1cm}

\noindent Let $U$ be a multiply connected domain of the Riemann sphere $\hC$ whose complement $\hC\setminus U$ has $N<\infty$ components. We show that every conformal map on $U$ can be written as a composition of $N$ maps conformal on simply connected domains. This improves on a recent result of D.E. Marshall, but our proof uses different ideas, and involves the uniformisation theorem for Riemann surfaces.

\section{Introduction}

It is well-known (see e.g.~\cite{Rudin}) that any function that is holomorphic on a multiply connected domain of the Riemann sphere $\hC$ whose complement has $N$ components (a $N$-tuply connected domain) can be represented as a sum of $N$ functions holomorphic on simply connected domains\footnote{Here and below, a domain is a nonempty connected open subset of the Riemann sphere.}. That is, if $A_i,\,i=1,\ldots,N$ are simply connected domains with pairwise disjoint complements, then a holomorphic function $f$ on $\cup_i A_i$ can be represented on $\cup_i A_i$ as $f=\sum_i f_i$ where each $f_i$ has a holomorphic continuation to $A_i$. The purpose of the present letter is to establish a similar fact for conformal maps (bijective conformal functions), where the sum is replaced by a factorisation into a composition of $N<\infty$ maps. In the case $N=2$, a related factorisation problem for increasing homeomorphisms of $\R$ (or of any other Jordan curve in $\hC$) is central in the context of conformal welding (see, e.g., \cite[p. 100]{Lehto}). The problem posed here also appeared, again in the case $N=2$, in a physically motivated study by the author of a differential structure on spaces of conformal maps \cite{D10}, where integral equations for the factors were found. In this context, holomorphic functions determine directions of small displacements of conformal maps, and the factorisation problem is the ``global'' equivalent of the ``local'' fact stated above about holomorphic functions. The $N>2$ version of the conformal welding problem was studied recently by D.E. Marshall in \cite{M11}, where Theorem 1 contains a result similar to ours: the conformal map of a domain bounded by $N$ quasicircles to a domain bounded by $N$ circles can be written as a composition of $N$ conformal maps of simply connected domains each of which extends to quasiconformal maps of the plane. Marshall points out that an arbitrary $N$-tuply connected domain can be mapped to a domain bounded by $N$ quasicircles using an additional $N$ maps of simply connected domains, resulting in a total of $2N$ maps in this case. Our theorem improves on the latter result by specifying that only $N$ factors are necessary and by lifting the condition that the image domain be bounded by circles; our proof makes use of different ideas.

We will show the following statement:
\begin{theorem}\label{theogenfact}
If $G$ is a conformal map on a $N$-tuply connected domain $U$ with $N<\infty$, then it is factorisable on this domain into a composition of $N$ factors $g_i,\;i=1,\ldots,N$ that are maps conformal on simply connected domains:
\beq\label{eqgenfact}
	G = g_1\circ g_2\circ\cdots\circ g_N \qquad \mbox{(on $U$)}.
\eeq
\end{theorem}
The main tool in the proof is Koebe's uniformisation theorem for Riemann surfaces. In particular, we will use the theorem that every simply connected compact Riemann surface is conformally equivalent to the Riemann sphere $\hC$. The idea relating the uniformisation theorem to the welding problem, which forms the basis of our proof, is briefly but nicely explained in a paper by Sharon and Mumford \cite{SM04}.

{\bf Acknowledgment.} I am very grateful to D.E. Marshall for constructive criticisms and encouragements to publish this result. I acknowledge support from the EPSRC First Grant ``From conformal loop ensembles to conformal field theory" EP/H051619/1.

\section{Proof of the theorem}

By definition, $D$ is a $n$-tuply connected domain if $\hC \setminus D = \cup_{j=1}^n E_j$ where each $E_j$ is closed and connected in $\hat C$, and where $E_j,\;j=1,\ldots,n$ are pairwise disjoint. Then\footnote{$D=\cap_{j=1}^n A_j$ is immediate, and $\p D=\cup_{j=1}^n \p A_j$ follows from this and the fact that $\hC\setminus A_j,\;j=1,\ldots,n$ are pairwise disjoint. Since $\p D$ has $n$ components \cite[V.14.5]{Newman}, each $\p A_j$ has exactly one component, hence each $A_j$ is simply connected \cite[VI.4.3]{Newman}. Note that a $n$-tuply connected domain according to our definition has connectivity $n-1$ according to \cite[VI.10.3]{Newman}.}, $A_j:=\hC\setminus E_j$ are simply connected domains with disjoint complements; $D=\cap_{j=1}^n A_j$; and $\p D = \cup_{j=1}^n \p A_j$.

For the proof of the theorem, we proceed by induction. Assume that there is an integer $n\leq N$, a map $g$ conformal on an $n$-tuply connected domain $D$, and, if $n<N$, maps $g_k$ for $k=n+1,\ldots,N$ conformal on simply connected domains, such that, on $U$, we have $G = g\circ g_{n+1}\circ g_{n+2}\circ \cdots\circ g_N$ if $n<N$ and $G=g$ if $n=N$. This is true for $n=N$, with $D=U$, by the assumption of the theorem. We will show that the assumptions of the induction imply that there is a map $\t{g}$ conformal on a $(n-1)$-tuply connected domain, and a map $g_n$ conformal on a simply connected domain, such that $G = \t{g}\circ g_n\circ g_{n+1}\circ\cdots\circ g_N$ on $U$.

Recall that $D = \cap_{i=1}^n A_i$ where $A_i$ are simply connected domains with pairwise disjoint complements. Since $g$ is a homeomorphism, $g(D)$ has the same connectivity as that of $D$ \cite[VI.12.1]{Newman}, hence likewise $\t{D}:= g(D) = \cap_{i=1}^n \t{A}_i$ for simply connected domains $\t{A}_i,\;i=1,\ldots,n$ with pairwise disjoint complements.

Since $g$ is proper, the cluster values of $g$ on $\p D$ lie in $\p \t{D}$. Since the complements of $A_i,\;i=1,\ldots,n$ are pairwise disjoint, for every $i,j$ with $i\neq j$, there is a simple closed polygon in $D$ separating $\hC\setminus A_i$ from $\hC\setminus A_j$ \cite[VI.3.3]{Newman}. Two points in $D$ that are separated by a simple closed polygon in $D$ map under $g$ to points that are separated by the image of the polygon in $\t{D}$, because $g$ is a homeomorphism (the image of the simple closed polygon is a Jordan curve hence separating \cite[V.10.2]{Newman}, and lies in $g(D)$). Therefore, for all $i,j$ with $i\neq j$, the cluster values of $g$ on $\p A_i$ and those on $\p A_j$ must lie on different components of $\p \t{D}$. Hence there is a bijective map $i\mapsto j$ such that all cluster values of $g$ on $\p A_i$ lie in $\p \t{A}_j$. We may choose the indexing such that this map is the identity. Clearly, this implies that for every $i$, the cluster values of $g^{-1}$ on $\p\t{A}_i$ lie in $\p A_i$. Hence, whenever $(z_\ell:\ell\in\Z)$ is such that $z_\ell$ and $g(z_\ell)$ converge,
\beq\label{bij}
	\lim_{\ell\to\infty} g(z_{\ell}) \in \p \t{A}_i\quad\Leftrightarrow \quad 	\lim_{\ell\to\infty} z_\ell\in \p A_i.
\eeq

Let $\t{A}:=\cap_{i=1}^{n-1} \t{A}_i$; this is a $(n-1)$-tuply connected domain (it is obtained from $\t{D}$ by ``filling-in'' the $n^{\rm th}$ hole).

1. We construct a Riemann surface $\cal C$ with two coordinate charts and transition map $g$.

Let ${\cal A}_c$ be a copy of $\t{A}$ and ${\cal B}_c$ be a copy of $A_n$, and denote by $\psi_A:\t{A}\onto {\cal A}_c$ and $\psi_B:A_n\onto {\cal B}_c$ the canonical copy bijections. Let ${\cal C}_c$ be the disjoint union of ${\cal A}_c$ and ${\cal B}_c$, and define a map ${\cal G}:\psi_B(D)\subset {\cal B}_c \onto \psi_A(\t{D})\subset{\cal A}_c$ by
\[
	{\cal G}:=\psi_A \circ g\circ \psi_B^{-1}.
\]
We define an equivalence relation $\sim$ on ${\cal C}_c$ by: $z\sim w$ if and only if either $z\in\psi_B(D)$ and ${\cal G}(z)=w$ or $w\in \psi_B(D)$ and ${\cal G}(w)=z$. Let ${\cal C}:={\cal C}_c/\sim$. Denote by $\iota:{\cal C}_c\onto {\cal C}$ the quotient map, and let ${\cal A} := \iota({\cal A}_c)$, ${\cal B}:=\iota({\cal B}_c)$ and ${\cal D}:=\iota(\psi_B(D))=\iota(\psi_A(\t{D}))$. Clearly, ${\cal D}= {\cal A}\cap {\cal B}$ and ${\cal C} = {\cal A}\cup {\cal B}$. Further, denote by $\iota_A:=\iota|_{{\cal A}_c}$ and $\iota_B:=\iota|_{{\cal B}_c}$ the restrictions of $\iota$ to ${\cal A}_c$ and ${\cal B}_c$. Note that $\iota_A:{\cal A}_c\onto {\cal A}$ and $\iota_B:{\cal B}_c\onto {\cal B}$ are bijections. Then,
\beq\label{eqtheogf}
	\iota_A\circ {\cal G}\circ \iota_B^{-1}=\id \quad\mbox{on}\quad {\cal D}.
\eeq
Put on ${\cal A}_c$ and ${\cal B}_c$ the topologies induced by those on $\t{A}$ and $A_n$, respectively, by the copy bijections. Then, by restricting these topologies, ${\cal G}$ is a homeomorphism from $\psi_B(D)$ onto $\psi_A(\t{D})$. Hence, $\iota_A$ and $\iota_B$ induce topologies on ${\cal A}$ and ${\cal B}$, respectively, which agree on ${\cal D}$, and we may put on $\cal C$ the topology generated by those on $\cal A$ and $\cal B$.

The topological space ${\cal C}$ is Hausdorff. This is shown as follows. Note that any $p\in{\cal A}$ has neighbourhoods lying entirely in ${\cal A}$: the images of those in $\t{A}$ under the map $\iota_A\circ \psi_A$; likewise, every $p\in{\cal B}$ has neighbourhoods lying entirely in ${\cal B}$: the images of those in $A_n$ under the map $\iota_B\circ\psi_B$. Then, it is clear that any two points in ${\cal A}$, or any two points in ${\cal B}$, have distinct neighbourhoods. If $p\in {\cal A}$ and all neighbourhoods of $p$ intersect ${\cal B}$, then this must hold true for those neighbourhoods that lie entirely in ${\cal A}$, hence $(\psi_A^{-1}\circ\iota_A^{-1})(p)\in\cl{\t{D}}$ whence $p\in \cl{{\cal D}}$. Let $p_1\in {\cal A}-{\cal B}$ and $p_2\in{\cal B}-{\cal A}$. If $p_1\not \in \p {\cal D}$, then the previous statement shows that $p_1$ has neighbourhoods not intersecting ${\cal B}$. Since there are neighbourhoods of $p_2$ lying entirely in ${\cal B}$, $p_1$ and $p_2$ have distinct neighbourhoods. A similar argument holds if instead $p_2\not\in \p{\cal D}$. Hence assume that $p_1,p_2\in\p {\cal D}$. If all neighbourhoods of $p_1$ and $p_2$ intersect, then this holds true for neighbourhoods ${\cal N}_1$ and ${\cal N}_2$ of $p_1$ and $p_2$, respectively, that lie entirely in ${\cal A}$ and ${\cal B}$, respectively, hence such that ${\cal N}_1\cap {\cal N}_2\subset {\cal D}$. Then there are sequences in ${\cal D}$ that converge to both $p_2$ and $p_1$. This implies that there are sequences in $D$ that converge to $(\psi_B^{-1}\circ\iota_B^{-1})(p_2)\in \p A_n$ and that map under $g$ to sequences in $\t{D}$ that converge to $(\psi_A^{-1}\circ\iota_A^{-1})(p_1)\in\p \t{A}$, which is in contradiction with (\ref{bij}). Hence, $p_1$ and $p_2$ have disjoint neighbourhoods.

Since in ${\cal C}$ every point has neighbourhoods homeomorphic to $\R^2$, the paragraph above implies that $\cal C$ is a topological surface. Let $\phi_A:=\psi_A^{-1}\circ\iota_A^{-1}$ and $\phi_B:=\psi_B^{-1}\circ \iota_B^{-1}$. These are coordinate charts for this surface. By inverting (\ref{eqtheogf}), we see that, on $D$,
\beq\label{factor}
	\phi_A\circ \phi_B^{-1} = g.
\eeq
Hence, the transition map is conformal, so that ${\cal C}$ is a Riemann surface. Note that ${\cal A}$ and ${\cal B}$ are sub-Riemann surfaces.

2. We show that ${\cal C}$ is simply connected.

Let $\gamma$ be a simple loop with base point $p\in{\cal B}-{\cal A}$. It will be sufficient to prove that every such loop in ${\cal C}$ is null-homotopic.

If $\gamma \subset{\cal  B}$, then $\gamma$ is null-homotopic in ${\cal B}$ because ${\cal B}$ is homeomorphic to $A_n$, which is simply connected. Hence, assume that $\gamma \cap {\cal A}-{\cal B}\neq\emptyset$. Note that since ${\cal A}-{\cal B} = {\cal A}\setminus {\cal D}$, we find $\phi_A({\cal A}-{\cal B}) = \phi_A({\cal A}\setminus {\cal D}) = \t{A}\setminus g(D) = \hC\setminus\t{A}_n$.

There exists a simply connected open neighbourhood $E$ of $\hC\setminus \t{A}_n$ bounded by a simple closed polygon, such that the closure of $E$ lies in $\t{A}$. Indeed, since $\hC\setminus \t{A}_i,\;i=1,\ldots,n$ are pairwise disjoint, the closed sets $\hC\setminus \t{A}_i,\;i=1,\ldots,n-1$ lie in $\t{A}_n$. Since $\t{A}_n$ is connected, there exists a path $\alpha$ whose trace intersects all of $\hC\setminus \t{A}_i,\;i=1,\ldots,n-1$. Hence $[\alpha]\cup \cup_{i=1}^{n-1}(\hC\setminus \t{A}_i)$ and $\hC\setminus \t{A}_n$ are disjoint nonempty connected closed sets, whence by \cite[VI.3.3]{Newman} there exists a simple closed polygon in $\t{D}$ separating them. This polygon bounds a simply connected domain $E$ that contains $\hC\setminus \t{A}_n$ and whose closure lies in $\t{A}$; that is, $\hC\setminus \t{A}_n\subset E$ and $\cl{E}\subset \t{A}$.

Let ${\cal E}:=\phi_A^{-1}(E)$. We have $\p {\cal E} = \phi_A^{-1}(\p E)$ because ${\cal C}$ is Hausdorff and $\phi_A$ is a homeomorphism. Hence, ${\cal A}-{\cal B}\subset {\cal E}$ and $\cl{{\cal E}}\subset \cal A$. Clearly $p\not\in {\cal A}$ whence $p\not\in\cl{{\cal E}}$. Then, $\gamma\setminus \p {\cal E}$ is a countable union of arcs in ${\cal C}$, each of which being a cross-cut of either ${\cal E}$ or ${\cal C}\setminus {\cal E}$. Each arc $\alpha$ which intersects ${\cal E}$ maps under $\phi_A$ to a cross-cut of $E$. This cross-cut is homotopic, in ${\cl E}$, to a subarc of $\partial E$, since $E$ is simply connected. Hence, thanks to the homeomorphism $\phi_A^{-1}$, each such arc $\alpha$ is homotopic in $\cl{{\cal E}}$ to a subarc of $\p {\cal E}$. This implies that $\gamma$ is homotopic to a curve $\widetilde \gamma$ which does not intersect $\cal A - \cal B$. But then $\widetilde \gamma \subset \cal B$, which is simply connected. This implies that $\widetilde \gamma$, and hence $\gamma$, is null-homotopic.

3. We show that ${\cal C}$ is compact.

We use the fact that the union of two compact subsets of a topological space is compact. Consider the subset ${\cal E}\subset {\cal C}$ constructed above. Since $\cl{E}\subset \t{A}$ is compact, $\phi_A$ is a homeomorphism and ${\cal C}$ is Hausdorff, then $\phi_A^{-1}(\cl{E})$ is compact. Hence $\cl{{\cal E}}\subset {\cal A}$ is compact, because $\phi_A^{-1}(\cl{E}) = \cl{{\cal E}}$. Similarly, $\hC\setminus E \subset \t{A}_n$ is compact and $\phi_B$ is a homeomorphism, whence $\phi_B^{-1}(\hC\setminus E) = {\cal C}\setminus {\cal E}\subset {\cal B}$ is compact. Clearly $\cl{{\cal E}} \cup ({\cal C}\setminus {\cal E}) = {\cal C}$, hence ${\cal C}$ is compact.

4. We complete the induction using the uniformisation theorem.

By the uniformisation theorem, there exists a conformal map $u:{\cal C} \onto \hC$. Hence, defining $\t{g} := \phi_A\circ u^{-1}$ and $g_n := u \circ \phi_B^{-1}$, Equation (\ref{factor}) implies $g = \t{g}\circ g_n$ on $D$, where $g_n$ is a conformal map of the simply connected domain $A_n$, and $\t{g}$ is a conformal map of the $(n-1)$-tuply connected domain $u({\cal A})$.

\end{document}